\documentclass{amsproc}%
\usepackage{graphicx}
\usepackage{amscd}
\usepackage{amsmath}
\usepackage{amsfonts}
\usepackage{amssymb}%
\theoremstyle{plain}

\numberwithin{equation}{section}

\begin{document}
\title[Subrings which are closed with respect to taking the inverse]{Subrings which are closed with respect to taking the inverse}
\author{Jen\H{o} Szigeti}
\address{Institute of Mathematics, University of Miskolc, Miskolc, Hungary 3515}
\email{jeno.szigeti@uni-miskolc.hu}
\author{Leon van Wyk}
\address{Department of Mathematical Sciences, Stellenbosch University\\
P/Bag X1, Matieland 7602, South Africa }
\email{lvw@sun.ac.za}
\thanks{\noindent The first author was supported by OTKA of Hungary No. T043034 and
K61007. The second author was supported by the National Research Foundation of
South Africa under Grant No. UID 61857. Any opinion, findings and conclusions
or recommendations expressed in this material are those of the authors and
therefore the National Research Foundation does not accept any liability in
regard thereto.}

\begin{abstract}
Let $S$ be a subring of the ring $R$. We investigate the question of whether
$S\cap U(R)=U(S)$ holds for the units. In many situations our answer is
positive. There is a special emphasis on the case when $R$ is a full matrix
ring and $S$ is a structural subring of $R$ defined by a reflexive and
transitive relation.

\end{abstract}
\maketitle

\noindent1. INTRODUCTION

\bigskip

\noindent Throughout the paper a ring $R$\ means a ring with identity, and all
subrings inherit the identity. The group of units in $R$ is denoted by $U(R)$
and the centre of $R$ is denoted by $Z(R)$.

\noindent In general, if $S$ is a subring of the ring $R$ and $x\in S$ is an
invertible element in $R$, then $x^{-1}$ need not be in $S$. The aim of this
paper is to investigate the question of whether $S\cap U(R)=U(S)$ holds for a
subring $S\subseteq R$. For a structural matrix subring of a full matrix ring
this question was raised by Johan Meyer. A similar problem for the additive
subgroups of a division ring was considered in [3].

\noindent In Section 2 first we impose certain chain conditions on $S$ or on
$R$ to derive that $S\cap U(R)=U(S)$. Then we combine the chain conditions
with the assumption that $R$ is a PI-ring. Using the prime ideals of $R$ we
formulate a reduction theorem providing $S\cap U(R)=U(S)$. In Section 2 we
also deal with the subrings of a full matrix ring (over a Noetherian or a PI-ring).

\noindent Section 3 is devoted to the study of the structural matrix subring
$M_{n}(\theta,R)$\ of the full matrix ring $M_{n}(R)$ defined by a reflexive
and transitive relation $\theta$ on the set $\{1,2,...,n\}$. First we
reformulate the general results of Section 2 to see that $M_{n}(\theta,R)\cap
U(M_{n}(R))=U(M_{n}(\theta,R))$ holds for various base rings $R$. Then we get
the same equality for PI-rings. Finally we prove that $M_{n}(\theta,R)$ is
closed with respect to taking the adjoint (note that the adjoint always
exists, not as the inverse).

\noindent In proving our statements we use some classical and one recent
theorem concerning PI-rings.

\noindent Section 4 contains an example (based on a classical construction of
Jacobson) indicating that the Noetherian and the PI conditions play an
adequate role in our development. Since any non-Dedekind-finite ring can
appear as a base ring in our example, we can use the results of Section 3 to
derive some more or less known statements about Dedekind-finite rings. The
authors are grateful to Peter P. P\'{a}lfy for his help in Section 4.

\bigskip

\noindent2. CHAIN AND PI\ CONDITIONS

\bigskip

\noindent It is known that a ring $R$ is called strongly $\pi$-regular if for
every $x\in R$ the DCC holds for the left ideals $Rx^{i}$, $i\geq1$.

\bigskip

\noindent\textbf{Proposition 2.1.}\textit{ Let }$R$\textit{ be an arbitrary
ring and let }$S$\textit{ be a strongly }$\pi$\textit{-regular subring of }%
$R$\textit{. If }$x\in S$\textit{ is invertible in }$R$\textit{, then }%
$x^{-1}\in S$\textit{.}

\bigskip

\noindent\textbf{Proof.} The DCC for the left ideals $Sx^{i}$, $i\geq1$, of
$S$\ gives that $Sx^{k}=Sx^{k+1}$ for some $k\geq1$. Thus%
\[
x^{k}=sx^{k+1}%
\]
for some $s\in S$, whence we obtain that $x^{-1}=s$ is in $S$. $\square$

\bigskip

\noindent\textbf{Proposition 2.2.}\textit{ Let }$R$\textit{ be a ring integral
over a central subring }$C\subseteq Z(R)$\textit{ and let }$C\subseteq
S\subseteq R$\textit{ be a subring. If }$x\in S$\textit{ is invertible in }%
$R$\textit{, then }$x^{-1}\in S$\textit{.}

\bigskip

\noindent\textbf{Proof.} The integrality gives that%
\[
x^{-k}+c_{k-1}x^{-(k-1)}+\cdots+c_{1}x^{-1}+c_{0}=0
\]
holds for $x^{-1}\in R$, where $k\geq1$ and $c_{k-1},...,c_{1},c_{0}\in C$.
Thus%
\[
x^{-1}=-(c_{k-1}+c_{k-2}x\cdots+c_{1}x^{k-2}+c_{0}x^{k-1})
\]
is in $S$. $\square$

\bigskip

\noindent\textbf{Proposition 2.3.}\textit{ Let }$S$\textit{ be a subring of
the ring }$R$\textit{ such that }$R$\textit{ is Noetherian as a left }%
$S$\textit{-module. If }$x\in S$\textit{ is invertible in }$R$\textit{, then
}$x^{-1}\in S$\textit{.}

\bigskip

\noindent\textbf{Proof.} The ACC for the $S$-submodules%
\[
H_{k}=\underset{i=1}{\overset{k}{\sum}}Sx^{-i},k\geq1,
\]
of the left $S$-module $_{S}R$ gives that%
\[
\underset{i=1}{\overset{k}{\sum}}Sx^{-i}=\underset{i=1}{\overset{k+1}{\sum}%
}Sx^{-i}%
\]
for some $k\geq1$. Thus%
\[
x^{-(k+1)}=s_{1}x^{-1}+s_{2}x^{-2}+\cdots+s_{k}x^{-k}%
\]
with $s_{1},s_{2},...,s_{k}\in S$, whence right multiplication by $x^{k}$
gives that%
\[
x^{-1}=s_{1}x^{k-1}+s_{2}x^{k-2}+\cdots+s_{k-1}x+s_{k}%
\]
is in $S$. $\square$

\bigskip

\noindent\textbf{Theorem 2.4.}\textit{ Let }$R$\textit{ be a prime PI-ring
such that }$Z(R)$\textit{ is Noetherian. If }$Z(R)\subseteq S\subseteq
R$\textit{ is a subring, then }$S\cap U(R)=U(S)$\textit{.}

\bigskip

\noindent\textbf{Proof.} A theorem of Formanek (see page 109 in Vol. II of
[4]) ensures that $R$ is a Noetherian $Z(R)$-module. The condition
$Z(R)\subseteq S$\ ensures that an $S$-submodule of the left $S$-module
$_{S}R$ is a $Z(R)$-submodule of $R$, whence we obtain that $R$ is Noetherian
as a left $S$-module. Thus Proposition 2.3 can be applied to the pair of rings
$S\subseteq R$. $\square$

\bigskip

\noindent\textbf{Theorem 2.5.}\textit{ Let }$R$\textit{ be a prime PI-ring
such that }$Z(R)$\textit{ is Noetherian. If }$S$\textit{ is a subring of
}$M_{n}(R)$\textit{\ such that }$\{rI\mid r\in Z(R)\}\subseteq S$\textit{,
then }$S\cap U(M_{n}(R))=U(S)$\textit{.}

\bigskip

\noindent\textbf{Proof.} Since $M_{n}(R)$ is also a prime PI-ring (see page
110 in Vol. II of [4]) with Noetherian centre%
\[
Z(M_{n}(R))=\{rI\mid r\in Z(R)\}\cong Z(R),
\]
the application of Theorem 2.4 gives the desired equality. $\square$

\bigskip

\noindent\textbf{Theorem 2.6.}\textit{ Let }$R$\textit{ be a left Noetherian
ring. If }$S$\textit{ is a subring of }$M_{n}(R)$\textit{\ such that
}$\{rI\mid r\in R\}\subseteq S$\textit{, then }$S\cap U(M_{n}(R))=U(S)$%
\textit{.}

\bigskip

\noindent\textbf{Proof.} Since $M_{n}(R)$ is a free left $R$-module (of rank
$n^{2}$), $M_{n}(R)$ is Noetherian as a left $R$-module. The condition
$\{rI\mid r\in R\}\subseteq S$\ ensures that an $S$-submodule of $_{S}%
M_{n}(R)$ is an $R$-submodule of $_{R}M_{n}(R)$, whence we obtain that
$M_{n}(R)$ is Noetherian as a left $S$-module. Thus Proposition 2.3 can be
applied to the pair of rings $S\subseteq M_{n}(R)$. $\square$

\bigskip

\noindent\textbf{Theorem 2.7.}\textit{ Let }$P_{i}\vartriangleleft R$\textit{,
}$1\leq i\leq t$\textit{, be a finite collection of ideals of the ring ~}%
$R$\textit{ such that the intersection }$P_{1}\cap P_{2}\cap\cdots\cap P_{t}%
$\textit{ is a nil ideal. For a subring }$S\subseteq R$\textit{\ consider the
subring }$S/P_{i}=\{s+P_{i}\mid s\in S\}\subseteq R/P_{i}$\textit{ of the
factor ring }$R/P_{i}$\textit{. If }$(S/P_{i})\cap U(R/P_{i})=U(S/P_{i}%
)$\textit{ for all }$i\in\{1,2,...,t\}$\textit{, then }$S\cap U(R)=U(S)$%
\textit{.}

\bigskip

\noindent\textbf{Proof.} Take an element $x\in S\cap U(R)$. Since $x+P_{i}\in
S/P_{i}$, our assumption gives that the inverse $(x+P_{i})^{-1}=x^{-1}+P_{i}$
is in $S/P_{i}$. Thus $x^{-1}+P_{i}=s_{i}+P_{i}$ for some $s_{i}\in S$. In
view of%
\[
1-xs_{i}=x(x^{-1}-s_{i})\in P_{i},
\]
we obtain that%
\[
(1-xs_{1})(1-xs_{2})\cdots(1-xs_{t})\in P_{1}P_{2}\cdots P_{t}\subseteq
P_{1}\cap P_{2}\cap\cdots\cap P_{t}%
\]
is nilpotent. Clearly,%
\[
(1-xs_{1})(1-xs_{2})\cdots(1-xs_{t})=1-xs
\]
for some $s\in S$, whence%
\[
0=(1-xs)^{k}=1+\tbinom{k}{1}(-xs)+\cdots+\tbinom{k}{k}(-xs)^{k}%
\]
follows for some integer $k\geq1$. Consequently%
\[
x^{-1}=\tbinom{k}{1}s-\tbinom{k}{2}s(xs)+\cdots+(-1)^{k+1}\tbinom{k}%
{k}s(xs)^{k-1}\in S.\text{ }\square
\]

\bigskip

\noindent\textbf{Theorem 2.8.}\textit{ Let }$R$\textit{ be a ring with ACC on
ideals and }$S\subseteq R$\textit{\ be a subring such that }$(\!S/P\!)\cap
U(\!R/P\!)\!=\!U(\!S/P\!)$\textit{ for all prime ideals }$P\vartriangleleft
R$\textit{. Then }$S\cap U(\!R\!)\!=\!U(\!S\!)$\textit{.}

\bigskip

\noindent\textbf{Proof.} The ACC ensures that the prime radical of $R$\ is a
finite intersection of prime ideals (see page 364 in Vol. I of [4]):%
\[
\text{rad}(R)=P_{1}\cap P_{2}\cap\cdots\cap P_{t}.
\]
Since rad$(R)$ is nil, we can apply Theorem 2.7 to get the desired equality.
$\square$

\bigskip

\noindent In the rest of this section we shall make use of a Lie nilpotent $R$
of index $m$\ as the underlying ring in $M_{n}(R)$, in other words a ring $R$
satisfying the identity%
\[
\lbrack\lbrack\lbrack...[[x_{1},x_{2}],x_{3}],...],x_{m}],x_{m+1}]=0,
\]
with $[x,y]=xy-yx$. The following theorem can easily be obtained from
Proposition 4.1 and Theorem 4.2 in [5]:

\bigskip

\noindent\textbf{Theorem 2.T.}\textit{ If }$R$\textit{ is a ring satisfying
the identity}%
\[
\lbrack\lbrack\lbrack...[[x_{1},x_{2}],x_{3}],...],x_{m}],x_{m+1}]=0
\]
\textit{and }$A\in M_{n}(R)$\textit{, then a left Cayley-Hamilton identity}%
\[
\lambda_{d}A^{d}+\lambda_{d-1}A^{d-1}+\cdots+\lambda_{1}A+\lambda_{0}I=0
\]
\textit{holds for }$A$\textit{, with }$d=n^{m}$\textit{, }$\lambda_{d}%
\in\mathbb{Z}\smallsetminus\{0\}$\textit{ and }$\lambda_{i}\in R$\textit{,
}$0\leq i\leq d$\textit{.}

\bigskip

\noindent\textbf{Theorem 2.9.}\textit{ Let }$R$\textit{ be a ring such that
}$\mathbb{Z\smallsetminus\{}0\mathbb{\}}\subseteq U(R)$\textit{ and }%
$R$\textit{ satisfies the identity}%
\[
\lbrack\lbrack\lbrack...[[x_{1},x_{2}],x_{3}],...],x_{m}],x_{m+1}]=0.
\]
\textit{If }$S$\textit{ is a subring of }$M_{n}(R)$\textit{\ such that
}$\{rI\mid r\in R\}\subseteq S$\textit{, then }$S\cap U(M_{n}(R))=U(S)$%
\textit{.}

\bigskip

\noindent\textbf{Proof.} If $A\in S\cap U(M_{n}(R))$, then Theorem 2.T
provides a left Cayley-Hamilton identity for $A^{-1}$\ of the form%
\[
\gamma_{d}A^{-d}+\gamma_{d-1}A^{-(d-1)}+\cdots+\gamma_{1}A^{-1}+\gamma
_{0}I=0.
\]
Since $\gamma_{d}\in\mathbb{Z}\smallsetminus\{0\}$ is in $U(R)$, right
multiplication by $A^{d-1}$ and then left multiplication by $\gamma_{d}^{-1}$
gives that%
\[
A^{-1}=-\gamma_{d}^{-1}(\gamma_{d-1}I+\gamma_{d-2}A+\cdots+\gamma_{1}%
A^{d-2}+\gamma_{0}A^{d-1})
\]
is in $S$. $\square$

\bigskip

\noindent\textbf{Corollary 2.10.}\textit{ Let }$R$\textit{ be a commutative
ring. If }$S$\textit{ is a subring of }$M_{n}(R)$\textit{\ such that
}$\{rI\mid r\in R\}\subseteq S$\textit{, then }$S\cap U(M_{n}(R))=U(S)$%
\textit{.}

\bigskip

\noindent\textbf{Proof.} We have $d=n$ and $\gamma_{n}=1$ in the classical
Cayley-Hamilton identity. $\square$

\bigskip

\noindent3. STRUCTURAL\ MATRIX RINGS

\bigskip

\noindent The class of structural matrix rings has been studied extensively,
see for example, [1] and [2]. For a reflexive and transitive binary relation
$\theta$ on the set $\{1,2,...,n\}$, the structural matrix subring
$M_{n}(\theta,R)$\ of the full matrix ring $M_{n}(R)$ is defined as follows:%
\[
M_{n}(\theta,R)=\{[a_{i,j}]\in M_{n}(R)\mid a_{i,j}=0\text{ if }%
(i,j)\notin\theta\}.
\]
Henceforth $\theta$ is a reflexive and transitive binary relation on
$\{1,2,...,n\}$. In the next three theorems we collect the consequences of
Theorems 2.6, 2.9 and Corollary 2.10.

\bigskip

\noindent\textbf{Theorem 3.1.}\textit{ If }$R$\textit{ is a left Noetherian
ring and }$A\in M_{n}(\theta,R)$\textit{ is invertible in }$M_{n}(R)$\textit{,
then }$A^{-1}\in M_{n}(\theta,R)$\textit{.}

\bigskip

\noindent\textbf{Theorem 3.2.} \textit{Let }$R$\textit{ be a Lie nilpotent
ring such that }$\mathbb{Z\smallsetminus\{}0\mathbb{\}}\subseteq
U(R)$\textit{. If }$A\in M_{n}(\theta,R)$\textit{ is invertible in }$M_{n}%
(R)$\textit{, then }$A^{-1}\in M_{n}(\theta,R)$\textit{.}

\bigskip

\noindent\textbf{Theorem 3.3.} \textit{If }$R$\textit{ is a commutative ring
and }$A\in M_{n}(\theta,R)$\textit{ is invertible in }$M_{n}(R)$\textit{, then
}$A^{-1}\in M_{n}(\theta,R)$\textit{.}

\bigskip

\noindent For arbitrary rings $R$ and $T$, let $\mu:R\longrightarrow T$ be a
ring homomorphism with $\mu(1)=1$ and consider the induced ring homomorphism%
\[
\mu_{n}:M_{n}(R)\longrightarrow M_{n}(T).
\]
Then the containment $\mu_{n}(M_{n}(\theta,R))\subseteq M_{n}(\theta,T)$ is
obvious. In addition, if $\mu$ is injective and $\mu_{n}(A)\in M_{n}%
(\theta,T)$, then $A\in M_{n}(\theta,R)$.

\noindent Let $\theta^{\prime}$ and $\theta^{\prime\prime}$ be reflexive and
transitive binary relations on the sets $\{1,2,...,n\}$ and $\{1,2,...,m\}$
respectively. Then it is evident that%
\[
M_{n}(\theta^{\prime},M_{m}(\theta^{\prime\prime},R))\cong M_{nm}%
(\overline{\theta},R)
\]
for every ring $R$, where $\overline{\theta}$ is the reflexive and transitive
binary relation on the set $\{1,2,...,nm\}$ defined by%
\[
\overline{\theta}=\{(i,j)\mid(\left\lceil i/m\right\rceil ,\left\lceil
j/m\right\rceil )\in\theta^{\prime}\text{ and }(i_{(m)},j_{(m)})\in
\theta^{\prime\prime}\},
\]
where $\left\lceil \cdot\right\rceil $ denotes the ceiling and%
\begin{align*}
i_{(m)}  &  \equiv i\text{ }(\operatorname{mod}m)\text{ and }1\leq i_{(m)}\leq
m,\\
j_{(m)}  &  \equiv j\text{ }(\operatorname{mod}m)\text{ and }1\leq j_{(m)}\leq
m.
\end{align*}
We are now in a position to state the following result.

\bigskip

\noindent\textbf{Theorem 3.4.}\textit{ Let }$C$\textit{ be a Noetherian
subring of }$Z(R)$\textit{ such that }$R$\textit{ is a PI-algebra over }%
$C$\textit{. If }$A\in M_{n}(\theta,R)$\textit{ and }$A$\textit{\ is
invertible in }$M_{n}(R)$\textit{, then }$A^{-1}$\textit{ is in }$M_{n}%
(\theta,R)$.

\bigskip

\noindent\textbf{Proof.} Let $A=[a_{i,j}]$, $A^{-1}=[b_{i,j}]$ and consider
the $C$-subalgebra%
\[
D=C\left\langle a_{i,j},b_{i,j}\mid1\leq i,j\leq n\right\rangle
\]
of $R$. Since $R$ is PI over $C$, the same holds for $D$. The theorem of
Razmyslov-Kemer-Braun (see page 151 in Vol. II of [4]) ensures that the upper nilradical

\noindent$N=$ Nil$(D)$ of the affine (finitely generated) $C$-algebra $D$ is
nilpotent: $N^{k}=\{0\}$ for some integer $k\geq1$. We know that $D/N$ can be
embedded in a full matrix ring $M_{m}(E)$ over a commutative ring $E$ (see
page 98 in Vol. II of [4]). Thus $\varphi=\mu\circ\varepsilon$ induces a ring
homomorphism%
\[
\varphi_{n}:M_{n}(D)\longrightarrow M_{n}(M_{m}(E)),
\]
where $\varepsilon:D\longrightarrow D/N$ is the natural surjection and
$\mu:D/N\longrightarrow M_{m}(E)$ is our embedding. Since $A$ is invertible in
$M_{n}(D)$, $\varphi_{n}(A)$ is invertible in $M_{n}(M_{m}(E))$. The argument
above gives that $M_{n}(\theta,M_{m}(E))\cong M_{nm}(\overline{\theta},E)$. By
assumption $A\in M_{n}(\theta,D)$, and so $\varphi_{n}(A)\in M_{n}%
(\theta,M_{m}(E))$ can be viewed as a matrix in $M_{nm}(\overline{\theta},E)$
having an inverse in $M_{nm}(E)$. Theorem 3.3 shows that the inverse
$(\varphi_{n}(A))^{-1}$, viewed as an $nm\times nm$ matrix over $E$, is in
$M_{nm}(\overline{\theta},E)$. As%
\[
\mu_{n}(\varepsilon_{n}(A^{-1}))=\varphi_{n}(A^{-1})=(\varphi_{n}(A))^{-1},
\]
we conclude that $\varepsilon_{n}(A^{-1})\in M_{n}(\theta,D/N)$ (see the above
observations preceding Theorem 3.4). Thus $\varepsilon(b_{i,j})=b_{i,j}+N=0$
holds in $D/N$ for all $(i,j)\notin\theta$. Define an $n\times n$ matrix
$W=[w_{i,j}]$ over $N$ as follows:%
\[
w_{i,j}=\left\{
\begin{array}
[c]{c}%
b_{i,j}\text{ if }(i,j)\notin\theta\\
0\text{ if }(i,j)\in\theta\text{ \ }%
\end{array}
\right.  .
\]
Now $A^{-1}-W\in M_{n}(\theta,D)$ and%
\[
I-AW=A(A^{-1}-W)\in M_{n}(\theta,D),
\]
whence $AW\in M_{n}(\theta,D)$ follows. Clearly, $W\in M_{n}(N)$ implies that
$AW\in M_{n}(N)$ and hence $(AW)^{k}=0$. In view of $A^{-1}-W=A^{-1}(I-AW)$,
we obtain that%
\[
A^{-1}=(A^{-1}-W)(I-AW)^{-1}=(A^{-1}-W)(I+(AW)+\cdots+(AW)^{k-1})
\]
is in $M_{n}(\theta,D)\subseteq M_{n}(\theta,R)$. $\square$

\bigskip

\noindent We note that the final calculations in the above proof can be
omitted by applying Theorem 2.7 to $M_{n}(\theta,D)\subseteq M_{n}(D)$ and
$P_{1}=M_{n}(N)\vartriangleleft M_{n}(D)$.

\bigskip

\noindent\textbf{Corollary 3.5.}\textit{ Let }$R$\textit{ be a PI-ring (a
PI-algebra over }$\mathbb{Z\subseteq}Z(R)$\textit{). If }$A\in~M_{n}%
(\theta,R)$\textit{ and }$A$\textit{\ is invertible in }$M_{n}(R)$\textit{,
then }$A^{-1}$\textit{ is in }$M_{n}(\theta,R)$.

\bigskip

\noindent Recall that in case $R$ is commutative, then a matrix $A=[a_{i,j}]$
$\in$ $M_{n}(R)$ is invertible if and only if $\det(A)$ $\in$ $U(R)$, in which
case%
\[
A^{-1}=(\det(A))^{-1}\text{adj}(A).
\]
For the classical adjoint matrix adj$(A)=[b_{r,s}]$ we have%
\[
b_{r,s}=\underset{\rho\in\text{S}_{n},\rho(s)=r}{\sum}\text{sgn}%
(\rho)a_{1,\rho(1)}\cdots a_{s-1,\rho(s-1)}a_{s+1,\rho(s+1)}\cdots
a_{n,\rho(n)},
\]
where the sum is taken over all permutations $\rho$\ of the set
$\{1,2,...,n\}$ with $\rho(s)=r$.

\noindent If $R$ is an arbitrary ring (not necessarily commutative), then the
preadjoint $A^{\ast}=[a_{r,s}^{\ast}]\in M_{n}(R)$ of $A=[a_{i,j}]\in
M_{n}(R)$ was defined as follows in [5]:%
\[
a_{r,s}^{\ast}=\underset{\tau,\rho}{\sum}\text{sgn}(\rho)a_{\tau(1),\rho
(\tau(1))}\cdots a_{\tau(s-1),\rho(\tau(s-1))}a_{\tau(s+1),\rho(\tau
(s+1))}\cdots a_{\tau(n),\rho(\tau(n))},
\]
where the sum is taken over all permutations $\tau$ of the set
$\{1,...,s-1,s+1,...,n\}$ and all permutations $\rho$\ of the set
$\{1,2,...,n\}$ with $\rho(s)=r$. If $R$ is commutative, then $A^{\ast
}=(n-1)!$adj$(A)$.

\bigskip

\noindent\textbf{Theorem 3.6.} \textit{If }$R$\textit{ is an arbitrary ring
and }$A\in M_{n}(\theta,R)$\textit{, then }$A^{\ast}\in M_{n}(\theta
,R)$\textit{.}

\bigskip

\noindent\textbf{Proof.} Let $1\leq r,s\leq n$, with $(r,s)\notin\theta$. We
prove that $a_{r,s}^{\ast}=0$. Take a permutation $\tau$ of the set
$\{1,...,s-1,s+1,...,n\}$ and a permutation $\rho$\ of the set $\{1,2,...,n\}$
with $\rho(s)=r$.

\noindent We claim that $(\tau(i),\rho(\tau(i)))\notin\theta$ for some
$i\in\{1,...,s-1,s+1,...,n\}$. Suppose the contrary, that is $(j,\rho
(j))\in\theta$ for all $j\in\{1,...,s-1,s+1,...,n\}$. Consider the cycle%
\[
(r,\rho(r),...,\rho^{t}(r))
\]
of the permutation $\rho$ (of length $t+1$ say). Since $\rho(s)=r$, it follows
that $\rho^{t}(r)=s$. The reflexivity of $\theta$ ensures that $r\neq s$, and
so%
\[
(r,\rho(r)),(\rho(r),\rho^{2}(r)),...,(\rho^{t-1}(r),s))\in\theta.
\]
The transitivity of $\theta$ implies that $(r,s)\in\theta$; a contradiction.
Thus $a_{\tau(i),\rho(\tau(i))}=0$ for some $i\in\{1,...,s-1,s+1,...,n\}$.
Consequently, each product%
\[
a_{\tau(1),\rho(\tau(1))}\cdots a_{\tau(s-1),\rho(\tau(s-1))}a_{\tau
(s+1),\rho(\tau(s+1))}\cdots a_{\tau(n),\rho(\tau(n))}%
\]
in the summation for $a_{r,s}^{\ast}$\ is zero, whence we obtain that
$a_{r,s}^{\ast}=0$. $\square$

\bigskip

\noindent\textbf{Corollary 3.7.}\textit{ If }$R$\textit{ is a commutative ring
and }$A\in M_{n}(\theta,R)$\textit{, then }adj$(A)\in M_{n}(\theta
,R)$\textit{.}

\bigskip

\noindent\textbf{Proof.} Comparing the definitions of adj$(A)$ and $A^{\ast}$
and ignoring $\tau$ in the proof of Theorem 3.6 shows that adj$(A)\in
M_{n}(\theta,R)$. $\square$

\bigskip

\noindent4. DEDEKIND-FINITE\ RINGS

\bigskip

\noindent A ring $R$ is called Dedekind-finite if $xy=1$ implies $yx=1$ for
all $x,y\in R$. The ring of linear transformations Hom$_{K}(V,V)$ of a vector
space $V$\ (over a field $K$) with a countably infinite basis $\{b_{1}%
,b_{2},...,b_{n},...\}\subseteq V$ is not Dedekind-finite. Define the linear
transformations $\alpha,\beta:V\longrightarrow V$ on the elements of the given
basis as $\alpha(b_{i})=b_{i-1}$ for $i\geq2$, $\alpha(b_{1})=0$ and
$\beta(b_{i})=b_{i+1}$ for $i\geq1$, then $\alpha\beta=1$ and $\beta\alpha
\neq1$. Note that the ring Hom$_{K}(V,V)$ is not left (right) Noetherian and
not PI.

\noindent The following example shows that we can not drop the left (or right)
Noetherian condition in Theorem 3.1 and the PI condition in Theorem 3.4.

\bigskip

\noindent\textbf{Example 4.1.} Let $R$ be an arbitrary non-Dedekind-finite
ring with elements $x,y\in R$ such that $xy=1$ and $yx\neq1$. The inverse of
the upper triangular $2\times2$ matrix
\[
A=\left[
\begin{array}
[c]{cc}%
y & 1-yx\\
0 & x
\end{array}
\right]
\]
over $R$ is the lower triangular $2\times2$ matrix%
\[
A^{-1}=\left[
\begin{array}
[c]{cc}%
x & 0\\
1-yx & y
\end{array}
\right]  .
\]
Thus $A\in M_{2}(\theta,R)$ and $A^{-1}\notin M_{2}(\theta,R)$, where
$\theta=\{(1,1),(1,2),(2,2)\}$. $\square$

\bigskip

\noindent In view of Theorems 3.1 and 3.4 the following corollaries can easily
be obtained.

\bigskip

\noindent\textbf{Corollary 4.2.}\textit{ If }$R$\textit{ is a left Noetherian
ring, then }$R$\textit{ is Dedekind-finite.}

\bigskip

\noindent\textbf{Corollary 4.3.}\textit{ If }$C$\textit{ is a Noetherian
subring of }$Z(R)$\textit{ such that }$R$\textit{ is a PI-algebra over }%
$C$\textit{, then }$R$\textit{ is Dedekind-finite.}

\bigskip

\noindent REFERENCES

\bigskip

\begin{enumerate}
\item G. Abrams, J. Haefner, A. del R\'{\i}o,\textit{ The isomorphism problem
for incidence rings,} Pacific J. Math. Vol. 187, No. 2 (1999), 201-214.

\item S. D\u{a}sc\u{a}lescu, L. van Wyk,\textit{ Do isomorphic structural
matrix rings have isomorphic graphs?, }Proc. Amer. Math. Soc. Vol. 124, No. 5
(1996), 1385-1391.

\item D. Goldstein, R. Guralnick, L. Small, E. Zelmanov, \textit{Inversion
invariant additive subgroups of division rings,} Pacific J. Math. Vol. 227,
No. 2 (2006), 287-294.

\item L.H. Rowen, \textit{Ring Theory,} Vol. I, II (1988) Academic Press Inc.,
San\ Diego, London.

\item J. Szigeti,\textit{\ New determinants and the Cayley-Hamilton theorem
for matrices over Lie nilpotent rings,} Proc. Amer. Math. Soc. Vol. 125, No. 8
(1997), 2245-2254.
\end{enumerate}

\end{document}